\documentclass[12pt]{amsart}
\usepackage{amsmath,amsthm,amsfonts,amssymb,mathrsfs}
\date{\today}

\usepackage{hyperref}

\newtheorem{theorem}{Theorem}[section]

\newtheorem{proposition}[theorem]{Proposition}
\newtheorem{corollary}[theorem]{Corollary}

\newtheorem{lemma}[theorem]{Lemma}

\theoremstyle{definition}
\newtheorem{example}[theorem]{Example}
\newtheorem{remark}[theorem]{Remark}

\begin{document}

\title[On chains in $H$-closed topological pospaces]
{On chains in $H$-closed topological pospaces}
\author{Oleg~Gutik}
\address{Department of Mechanics and Mathematics,
Ivan Franko Lviv National University,
Universytetska 1, Lviv, 79000, Ukraine}
\email{o\_\,gutik@franko.lviv.ua, ovgutik@yahoo.com}

\author{Du\v{s}an~Pagon}
\address{Institute of Mathematics, Physics and Mechanics, and
Faculty of Natural Sciences and Mathematics, University of
Maribor, Jadranska 19, Ljubljana, 1000, Slovenia}
\email{dusan.pagon@uni-mb.si}
\author{Du\v{s}an~Repov\v{s}}
\address{Institute of Mathematics, Physics and Mechanics, and
Faculty of Mathematics and Physics University of Ljubljana,
P.~O.~B. Jadranska 19, Ljubljana, 1000, Slovenia}
\email{dusan.repovs@guest.arnes.si}



\begin{abstract}
We study chains in an $H$-closed topological partially ordered
space. We give sufficient conditions for a maximal chain $L$ in an
$H$-closed topological partially ordered space such that $L$
contains a maximal (minimal) element. Also we give sufficient
conditions for a linearly ordered topological partially ordered
space to be $H$-closed. We prove that any $H$-closed topological
semilattice contains a zero. We show that a linearly ordered
$H$-closed topological semilattice is an $H$-closed topological
pospace and show that in the general case this is not true. We
construct an example an $H$-closed topological pospace with a
non-$H$-closed maximal chain and give sufficient conditions that a
maximal chain of an $H$-closed topological pospace is an
$H$-closed topological pospace.
\end{abstract}

\maketitle


\section{Introduction}

In this paper all topological spaces will be assumed to be
Hausdorff. We shall follow the terminology of~\cite{CHK, CP,
Engelking1989, GHKLMS, Nachbin1965, Ward1954}. If $A$ is a subset
of a topological space $X$, then we denote the closure of the set
$A$ in $X$ by $\operatorname{cl}_X(A)$.
By a \emph{partial order} on a set $X$ we mean a reflexive,
transitive and anti-symmetric binary relation $\leqslant$ on $X$.
If the partial order  $\leqslant$ on a set $X$ satisfies the
following linearity law
\begin{equation*}
    \mbox{if } \; x,y\in X, \; \mbox{ then } \; x\leqslant y \;
    \mbox{ or } \; y\leqslant x,
\end{equation*}
then it is said to be a \emph{linear order}. We write $x<y$ if
$x\leqslant y$ and $x\neq y$, $x\geqslant y$ if $y\leqslant x$,
and $x\nleqslant y$ if the relation $\leqslant y$ is false.
Obviously if $\leqslant$ is a partial order or a linear order on a
set $X$ then so is $\geqslant$. A set endowed with a partial order
(resp., linear order) is called a \emph{partially ordered} (resp.,
\emph{linearly ordered}) set. If $\leqslant$ is a partial order on
$X$ and $A$ is a subset of $X$ then we denote
\begin{equation*}
    {\downarrow}A=\, \{ y\in X\mid y\leqslant x \; \mbox{ for some
    }\; x\in A\} \quad
\mbox{and} \quad
    {\uparrow}A=\, \{ y\in X\mid x\leqslant y \; \mbox{ for some
    }\; x\in A\}.
\end{equation*}
Also for any elements $a, b$ of a partial ordered set $X$ such
that $a\leqslant b$ we denote ${\uparrow}a={\uparrow}\{ a\}$,
${\downarrow}a={\downarrow}\{ a\}$, $[a,b]={\uparrow}
a\cap{\downarrow}b$ and $[a,b)=[a,b]\setminus\{ b\}$. A subset $A$
of a partially ordered set $X$ is called \emph{increasing}
(\emph{decreasing}) if $A={\uparrow}A$ ($A={\downarrow}A$).

A partial order $\leqslant$ on a topological space $X$ is said to
be \emph{lower} (\emph{upper}) \emph{semicontinuous} provided,
that whenever $x\nleqslant y$ ($y\nleqslant x$) in $X$, then there
exists an open set $U\ni x$ such that if $a\in U$ then
$a\nleqslant y$ ($y\nleqslant a$). A partial order is called
\emph{semicontinuous} if it is both upper and lower
semicontinuous. Next, it is said to be \emph{continuous} or
\emph{closed} provided, that whenever $x\nleqslant y$ in $X$,
there exist open sets $U\ni x$ and $V\ni y$ such that if $a\in U$
and $b\in V$ then $a\nleqslant b$. Clearly, the statement that the
partial order $\leqslant$ on $X$ is semicontinuous is equivalent
to the assertion that ${\uparrow}a$ and ${\downarrow}a$ are closed
subsets of $X$ for each $a\in X$. A topological space equipped
with a continuous partial order is called a \emph{topological
partially ordered space} or shortly \emph{topological pospace}. A
partial order $\leqslant$ on a topological space $X$ is continuous
if and only if the graph of $\leqslant$ is a closed subset in
$X\times X$~\cite[Lemma~1]{Ward1954}. Also a semicontinuous linear
order on a topological space is
continuous~\cite[Lemma~3]{Ward1954}.

A \emph{chain} of a partially ordered set $X$ is a subset of $X$
which is linearly ordered with respect to the partial order. A
\emph{maximal chain} is a chain which is properly contained in no
other chain. The Axiom of Choice  implies the existence of maximal
chains in any partially ordered set. Every maximal chain in a
topological pospace is a closed set~\cite[Lemma~4]{Ward1954}.

An element $y$ in a partially ordered set $X$ is called
\emph{minimal} (resp. \emph{maximal}) in $X$ whenever $x\leqslant
y$ (resp. $y\leqslant x$) in $X$ implies $y\leqslant x$ (resp.
$x\leqslant y$). Let $X$ and $Y$ be partially ordered sets. A map
$f\colon X\rightarrow Y$ is called \emph{monotone} (or
\emph{partially order preserving}) if $x\leqslant y$ implies
$f(x)\leqslant f(y)$ for any $x, y\in X$.

A Hausdorff topological space $X$ is called \emph{$H$-closed} if
$X$ is a closed subspace of every Hausdorff space in which it is
contained~\cite{AlexandroffUrysohn1923, AlexandroffUrysohn1929}. A
Hausdorff pospace $X$ is called \emph{$H$-closed} if $X$ is a
closed subspace of every Hausdorff pospace in which it is
contained. Obviously that the notion of $H$-closedness is a
generalization of compactness. For any element $x$ of a compact
topological pospace $X$ there exists a minimal element $y\in X$
and a maximal element $z\in X$ such that $y\leqslant x\leqslant
z$~\cite{GHKLMS}. Every chain in a compact topological pospace is
a compact subset and hence it contains minimal and maximal
elements. Also for any point $x$ of a compact topological pospace
$X$ there exists a base at $x$ which consists of open order-convex
subsets~\cite{Nachbin1965} (A non-empty set $A$ of a partially
ordered set is called \emph{order-convex} if $A$ is an
intersection of increasing and decreasing subsets). We are
interested in the following question: \emph{Under which condition
an $H$-closed topological pospace has similar properties as a
compact topological pospace?}

In this paper we establish chains in $H$-closed topological
pospace. We give sufficient conditions on a maximal chain $L$ in
an $H$-closed topological pospace such that $L$ contains the
maximal (minimal) element. Also we give sufficient conditions on a
linearly ordered topological pospace to be $H$-closed. We prove
that every $H$-closed topological semilattice contains a zero. We
show that a linearly ordered $H$-closed topological semilattice is
an $H$-closed topological pospace and show that in the general
case it is not true. We construct an example an $H$-closed
topological pospace with a non-$H$-closed maximal chain and give
sufficient conditions that a maximal chain of an $H$-closed
topological pospace is an $H$-closed topological pospace.


\section{On maximal and minimal elements of maximal chains in
$H$-closed topological pospaces}

A chain $L$ of a partially ordered set $X$ is called
\emph{down-directed} (resp. \emph{up-directed}) in $X$ if for any
$x\in X$ there exists $l\in L$ such that $l\leqslant x$ (resp.
$x\leqslant l$).

\begin{theorem}\label{theorem1.1} Any down-directed
chain $L$ of an $H$-closed topological pospace contains a minimal
element of $L$.
\end{theorem}

\begin{proof} Suppose there exists an $H$-closed topological
pospace $X$ with a down-directed chain $L$ such that $L$ does not
contain a minimal element in $L$.

Let $x\notin X$. We extend the partial order $\leqslant$ from $X$
onto $X^\ast=X\cup\{ x\}$ as follows:
\begin{enumerate}
    \item[1)] $x\leqslant x$;
    \item[2)] $x\leqslant y$ for each $y\in L$; and
    \item[3)] $x\leqslant z$ for $z\in X\setminus L$ if and only if
    there exists $y\in L$ such that $y\leqslant z$.
\end{enumerate}

Then by Lemma~1~\cite{Ward1954} for any $a, b\in X$ such that
$a\nleqslant b$ there exist open neighbourhoods $U(a)$ and $U(b)$
of the points $a$ and $b$ respectively, such that
\begin{equation*}
    U(a)={\uparrow}U(a), \quad U(b)={\downarrow}U(b), \quad
    \mbox{and} \quad U(a)\cap U(b)=\varnothing.
\end{equation*}
We define the family
\begin{equation*}
    \mathscr{U}=\{(U(a),U(b))\mid a\nleqslant b,\; a\in X,\; b\in L\}
\end{equation*}
as follows: $U(a)$ is an open neighbourhood of $a$ and $U(b)$ is
an open neighbourhood of $b$ such that
\begin{equation*}
    U(a)={\uparrow}U(a), \quad U(b)={\downarrow}U(b), \quad
    \mbox{and} \quad U(a)\cap U(b)=\varnothing.
\end{equation*}

We denote the topology on $X$ by $\tau$. On $X^\ast$ we determine
a topology $\tau^\ast$ as follows. For any point $y\in X$ the
bases of topologies $\tau^\ast$ and $\tau$ at the point $y$
coincide. For any $y\in X$ by $\mathscr{B}(y)$ we denote the base
of the topology $\tau$ at the point $y$. We put
\begin{equation*}
\mathscr{P}(x)=\big\{\{ x\}\cup U(b)\mid
(U(a),U(b))\in\mathscr{U},\; a\nleqslant b,\; a\in X,\; b\in
L\big\}
\end{equation*}
and
\begin{equation*}
\mathscr{B}(x)=\{ U_1\cap\cdots\cap U_n\mid
U_1,\ldots,U_n\in\mathscr{P}(x),\; n=1,2,3,\ldots\}.
\end{equation*}
Obviously, the conditions (BP1)--(BP3) of~\cite{Engelking1989}
hold for the family $\{{\mathscr B}(y)\}_{y\in X^*}$ and hence
${\mathscr B}(x)$ is a base of a topology $\tau^\ast$ at the point
$x$.

Further we shall show that $(X^\ast,\tau^\ast,\leqslant)$ is a
topological pospace. Let $y\in X$. Then $y\notin{\downarrow}x$. We
consider two cases $y\in L$ and $y\in X\setminus L$. In the first
case we have $x<y$, and since $L$ does not contain a minimal
element, there exists $b\in L$ such that $b<y$. We put $W(x)=\{
x\}\cup U(b)$, where $(U(y),U(b))\in\mathscr{U}$. Obviously,
\begin{equation*}
W(x)={\downarrow}W(x), \quad U(y)={\uparrow}U(y), \quad \mbox{and}
\quad W(x)\cap U(y)=\varnothing.
\end{equation*}

Let $y\in X\setminus L$. Then there exists $b\in L$ such that
$y\nleqslant b$. In the other case we have $y\leqslant a$ for all
$a\in L$, a contradiction to the fact that $L$ is a down-directed
chain in $X$. Then we put $W(x)=\{ x\}\cup U(b)$, where
$(U(y),U(b))\in\mathscr{U}$. Therefore we have
\begin{equation*}
W(x)={\downarrow}W(x), \quad U(y)={\uparrow}U(y), \quad \mbox{and}
\quad W(x)\cap U(y)=\varnothing.
\end{equation*}

Thus $(X^\ast,\tau^\ast,\leqslant)$ is a Hausdorff topological
pospace which contains $X$ as a dense subspace, a contradiction.
The obtained contradiction implies that $L$ contains a minimal
element.
\end{proof}

The proof of the following theorem is similar.

\begin{theorem}\label{theorem1.1a} Any up-directed
chain of an $H$-closed topological pospace contains a maximal
element.
\end{theorem}

A subset $F$ of topological pospace $X$ is said to be \emph{upper}
(resp. \emph{lower}) \emph{separated} if and only if for each
$a\in X\setminus{\uparrow}F$ (resp. $a\in
X\setminus{\downarrow}F$) there exist disjoint open neighbourhoods
$U$ of $a$ and $V$ of $F$ such that $U$ is decreasing (resp.
increasing) and $V$ is increasing (resp. decreasing) in $X$.

\begin{theorem}\label{theorem1.1b} Any maximal upper separated
chain $L$ of an $H$-closed topological pospace $X$ contains a
minimal element of $L$.
\end{theorem}

\begin{proof} Suppose to the contrary that  there exists an $H$-closed
topological pospace $X$ with a maximal upper separated chain $L$
such that $L$ does not contain a minimal element.

Let $x\notin X$. We extend the partial order $\leqslant$ from $X$
onto $X^\ast=X\cup\{ x\}$ as follows:
\begin{enumerate}
    \item[1)] $x\leqslant x$;
    \item[2)] $x\leqslant y$ for each $y\in L$; and
    \item[3)] $x\leqslant z$ for $z\in X\setminus L$ if and only if
    there exists $y\in L$ such that $y\leqslant z$.
\end{enumerate}

Then by Lemma~1~\cite{Ward1954} we can to define the family
\begin{equation*}
    \mathscr{U}=\{(U(a),U(b))\mid a\nleqslant b,\; a\in X,\; b\in L\}
\end{equation*}
as follows: $U(a)$ is an open neighbourhood of $a$ and $U(b)$ is
an open neighbourhood of $b$ such that
\begin{equation*}
    U(a)={\uparrow}U(a), \quad U(b)={\downarrow}U(b), \quad
    \mbox{and} \quad U(a)\cap U(b)=\varnothing.
\end{equation*}

Since $L$ is a upper separated chain for any $a\in
X\setminus{\uparrow}L$ such that $a\nleqslant l$ for each $l\in
L$, there exist an open neighbourhood $V_a(L)$ of $L$ and an open
neighbourhood $V(a)$ of $a$ such that
\begin{equation*}
    V_a(L)={\uparrow}V_a(L), \quad V(a)={\downarrow}V(a), \quad
    \mbox{and} \quad V_a(L)\cap V(a)=\varnothing.
\end{equation*}
We define the family
\begin{equation*}
    \mathscr{V}=\{(V_a(L),V(a))\mid a\nleqslant l \; \mbox{ for
    any }\; l\in L,\; a\in X\setminus{\uparrow}L\}
\end{equation*}
as follows: $V_a(L)$ is an open neighbourhood of the chain $L$ and
$V(a)$ is an open neighbourhood of the point  $a$ such that
\begin{equation*}
    V_a(L)={\uparrow}V_a(L), \quad V(a)={\downarrow}V(a), \quad
    \mbox{and} \quad V_a(L)\cap V(a)=\varnothing.
\end{equation*}

We denote the topology on $X$ by $\tau$. On $X^\ast$ we determine
a topology $\tau^\ast$ as follows. For any point $y\in X$ the
bases of topologies $\tau^\ast$ and $\tau$ at the point $y$
coincide. For any $y\in X$ by $\mathscr{B}(y)$ we denote the base
of the topology $\tau$ at the point $y$. We put
\begin{equation*}
\mathscr{U}(x)=\big\{\{ x\}\cup U(b)\mid
(U(a),U(b))\in\mathscr{U},\; a\nleqslant b,\; a\in X,\; b\in
L\big\},
\end{equation*}
\begin{equation*}
\mathscr{V}(x)=\big\{\{ x\}\cup V_a(L)\mid
(V_a(L),V(a))\in\mathscr{V},\; a\nleqslant l \; \mbox{ for any }\;
l\in L,\; a\in X\setminus{\uparrow}L\big\},
\end{equation*}
\begin{equation*}
    \mathscr{P}(x)=\mathscr{U}(x)\cup\mathscr{V}(x), \qquad
\mbox{and}
\end{equation*}
\begin{equation*}
\mathscr{B}(x)=\{ U_1\cap\cdots\cap U_n\mid
U_1,\ldots,U_n\in\mathscr{P}(x),\; n=1,2,3,\ldots\}.
\end{equation*}
Obviously, the conditions (BP1)--(BP3) of~\cite{Engelking1989}
hold for the family $\{{\mathscr B}(y)\}_{y\in X^\ast}$ and
hence\break $\{{\mathscr B}(y)\}_{y\in X^\ast}$ is a base of a
topology $\tau^\ast$ at the point $y\in X^\ast$. Since the chain
$L$ does not contain a minimal element, every finite intersection
of elements from the family $\mathscr{U}(x)$ contains infinitely
many points from the chain $L$, and since every set of the family
$\mathscr{V}(x)$ contains the chain $L$, we conclude that $x$ is
not an isolated point in $(X^\ast,\tau^\ast)$.

Further we shall show that $(X^\ast,\tau^\ast,\leqslant)$ is a
topological pospace. Let $y\in X$. We consider two cases $y\in L$
and $y\in X\setminus L$. In the first case we have $x<y$, and
since $L$ does not contains a minimal element there exists $b\in
L$ such that $b<y$. We put $W(x)=\{ x\}\cup U(b)$, where
$(U(y),U(b))\in\mathscr{U}(x)$. Obviously,
\begin{equation*}
W(x)={\downarrow}W(x), \quad U(y)={\uparrow}U(y), \quad \mbox{and}
\quad W(x)\cap U(y)=\varnothing.
\end{equation*}

Let $y\in X\setminus L$. If $y\nleqslant x$, then there exists
$b\in L$ such that $y\nleqslant b$. In other case we have
$y\leqslant a$ for all $a\in L$, a contradiction to the maximality
of the chain $L$. Then we put $W(x)=\{ x\}\cup U(b)$, where
$(U(y),U(b))\in\mathscr{U}(x)$. Therefore we have
\begin{equation*}
W(x)={\downarrow}W(x), \quad U(y)={\uparrow}U(y), \quad \mbox{and}
\quad W(x)\cap U(y)=\varnothing.
\end{equation*}

If $x\nleqslant y$ then the definition of the family
$\mathscr{V}(x)$ implies that there exist an open neighbourhood
$V(x)$ of the point $x$ in $X^\ast$ and an open neighbourhood
$V(y)$ of the point $y$ in $X^\ast$ such that
\begin{equation*}
V(x)={\uparrow}V(x), \quad V(y)={\downarrow}V(y), \quad \mbox{and}
\quad V(x)\cap V(y)=\varnothing.
\end{equation*}

Thus $(X^\ast,\tau^\ast,\leqslant)$ is a Hausdorff topological
pospace which contains $X$ as a dense subspace, a contradiction.
The obtained contradiction implies that $L$ contains a minimal
element.
\end{proof}

The proof of the following theorem is similar.

\begin{theorem}\label{theorem1.1c} Any maximal lower separated
chain $L$ of an $H$-closed topological pospace $X$ contains a
maximal element of $L$.
\end{theorem}

Similarly to \cite{McCartan1971, Priestley1972} we shall say that
a topological pospace $X$ is a \emph{$C_i$-space} (resp.
\emph{$C_d$-space}) if whenewer a subset $F$ of $X$ is closed,
${\uparrow}F$ (resp. ${\downarrow}F$) is a closed subset in $X$. A
maximal chain of a topological pospace $X$ is called an
\emph{$MCC_i$-chain} (resp. an \emph{$MCC_d$-chain}) in $X$ if
${\uparrow}L$ (resp. ${\downarrow}L$) is a closed subset in $X$.
Obviously, if a topological pospace $X$ is a $C_i$-space (resp.
$C_d$-space) then any maximal chain in $X$ is an $MCC_i$-chain
(resp. $MCC_d$-chain) in $X$.

A topological pospace $X$ is said to be \emph{upper} (resp.
\emph{lower}) \emph{regularly ordered} if and only if for each
closed increasing (resp. decreasing) subset $F$ in $X$ and each
element $a\notin F$, there exist disjoint open neighbourhoods $U$
of $a$ and $V$ of $F$ such that $U$ is decreasing (resp.
increasing) and $V$ is increasing (resp. decreasing) in
$X$~\cite{Green1968, ChoePark1979}. A topological pospace $X$ is
\emph{regularly ordered} if it is upper and lower regularly
ordered.

Theorem~\ref{theorem1.1b} implies

\begin{corollary}\label{corollary1.1b} Any maximal
$MCC_i$-chain of an $H$-closed upper regularly ordered topological
pospace contains a minimal element of $L$.
\end{corollary}

Theorem~\ref{theorem1.1c} implies

\begin{corollary}\label{corollary1.1c} Any maximal
$MCC_d$-chain of an $H$-closed lower regularly ordered topological
pospace contains a maximal element of $L$.
\end{corollary}

\section{Some remarks on $H$-closed topological semilattices}

A topological space $S$ that is algebraically a semigroup with a
continuous semigroup operation is called a {\em topological
semigroup}. A {\em semilattice} is a semigroup with a commutative
idempotent semigroup operation. A {\em topological semilattice} is
a topological semigroup which is algebraically a semilattice.

If $E$ is a semilattice, then the semilattice operation on $E$
determines the partial order $\leqslant$ on $E$:
\begin{equation*}
e\leqslant f \quad \text{ if and only if } \quad ef=fe=e.
\end{equation*}
This order is called {\em natural}. A semilattice $E$ is called
{\em linearly ordered} if the semilattice operation admits a
linear natural order on $E$. The natural order on a topological
semilattice $E$ admits the structure of topological pospace on $E$
(see: \cite[Proposition~VI-1.14]{GHKLMS}). Obviously, if $S$ is a
topological semilattice then ${\uparrow}e$ and ${\downarrow}e$ are
closed subsets in $S$ for any $e\in S$.

A topological semilattice $S$ is called \emph{$H$-closed} if it is
a closed subset in any topological semilattice which contains $S$
as a subsemilattice. Properties of $H$-closed topological
semilattices were established in \cite{ChuchmanGutik2007,
GutikRepovs, Stepp1975}.

\begin{theorem}\label{theorem3.1}
Every $H$-closed topological semilattice contains the smallest
idempotent.
\end{theorem}

\begin{proof}
Suppose to the contrary that there exists an $H$-closed
topological semilattice $E$ which does not contain the smallest
idempotent. Let $x\notin E$. We put $E^\ast=E\cup\{ x\}$, and
extend the semilattice operation from $E$ onto $E^\ast$ as
follows:
\begin{equation*}
xx=xe=ex=x \qquad \mbox{for all} \quad e\in E.
\end{equation*}

Since $E$ is a topological pospace, there exist by
Lemma~1~\cite{Ward1954} for any $a, b\in E$ such that $a\nleqslant
b$ open neighbourhoods $U(a)$ and $U(b)$ of the points $a$ and $b$
respectively such that
\begin{equation*}
    U(a)={\uparrow}U(a), \quad U(b)={\downarrow}U(b), \quad
    \mbox{and} \quad U(a)\cap U(b)=\varnothing.
\end{equation*}
We define the family
\begin{equation*}
    \mathscr{U}=\{(U(a),U(b))\mid a\nleqslant b,\; a, b\in X\}
\end{equation*}
as follows: $U(a)$ is an open neighbourhood of $a$ and $U(b)$ is
an open neighbourhood of $b$ such that
\begin{equation*}
    U(a)={\uparrow}U(a), \quad U(b)={\downarrow}U(b), \quad
    \mbox{and} \quad U(a)\cap U(b)=\varnothing.
\end{equation*}

We denote the topology on $E$ by $\tau$. On $E^\ast$ we determine
a topology $\tau^\ast$ as follows. For any point $y\in E$ the
bases of topologies $\tau^\ast$ and $\tau$ at the point $y$
coincide. For any $y\in E$ by $\mathscr{B}(y)$ we denote the base
of the topology $\tau$ at the point $y\in E$.  We put
\begin{equation*}
\mathscr{P}(x)=\{\{ x\}\cup U(b)\mid (U(a),U(b))\in\mathscr{U},\;
a\nleqslant b,\; a,b\in E\}
\end{equation*}
and
\begin{equation*}
\mathscr{B}(x)=\{ U_1\cap\cdots\cap U_n\mid
U_1,\ldots,U_n\in\mathscr{P}(x),\; n=1,2,3,\ldots\}.
\end{equation*}
Obviously, the conditions (BP1)--(BP3) of~\cite{Engelking1989}
hold for the family $\{{\mathscr B}(y)\}_{y\in E^\ast}$ and hence
${\mathscr B}(y)$ is a base of a topology $\tau^\ast$ at the point
$y\in E^\ast$.

We shall further show that $(E^\ast, \tau^\ast)$ is a topological
semilattice. Let $e$ be an arbitrary element of the semilattice
$E$ and let $U$ be an arbitrary open neighbourhood of $x$ in
$E^\ast$ such that $e\notin U$. Since the semilattice $E$ does not
contain the smallest idempotent, such open neighbourhood $U$ of
the point $x$ exists in $E^\ast$. Otherwise we put $V=U\cap U(b)$,
where $b\in U$, $b<e$, and $(U(e),U(b))\in\mathscr{U}$. Thus
$V={\downarrow}V$ and $e\notin V$. Then there exist
$a_1,\ldots,a_n,b_1,\ldots,b_n\in E$, $a_1\nleqslant b_1,\ldots,
a_n\nleqslant b_n$, $n=1,2,3,\ldots$, and open neighbourhoods
$U(a_1),\ldots,U(a_n), U(b_1),\ldots,U(b_n)$ of the points
$a_1,\ldots,a_n,b_1,\ldots,b_n$ in $E^\ast$, respectively, such
that
\begin{enumerate}
    \item[1)] $U(a_1)={\uparrow}U(a_1), \ldots,
    U(a_n)={\uparrow}U(a_n)$;
    \item[2)] $U(b_1)={\downarrow}U(b_1), \ldots,
    U(b_n)={\downarrow}U(b_n)$;
    \item[3)] $U(a_1)\cap U(b_1)=\varnothing, \ldots,U(a_n)\cap
    U(b_n)=\varnothing$; and
    \item[4)] $U=U(a_1)\cap \ldots U(a_n)$.
\end{enumerate}

Let $b_0\in U\setminus\{ x\}$ such that $b_0<e$. Such element
$b_0$ exists since the semilattice $E$ does not contain the
smallest idempotent. Since $E$ is a topological semilattice there
exist an open neighborhood $V(e)$ of the point $e$ and open
neighbourhoods $W(b_0)$ and $V(b_0)$ of the point $b_0$ such that
\begin{equation*}
V(e)V(b_0)\subseteq W(b_0)\subseteq U.
\end{equation*}
Then there exist an open neighbourhood $U(e)$ of the point $e$ and
open neighbourhood $U(b_0)$ of the point $b_0$ such that
$(U(e),U(b_0))\in\mathscr{U}$. We put $O(e)=U(e)\cap V(e)$ and
$W=U\cap U(b_0)$. Since $E$ is a topological semilattice and
$W={\downarrow}W$, we have
\begin{equation*}
O(e)W\subseteq W\subseteq U.
\end{equation*}
Obviously $W\in\mathscr{B}(x)$.

Therefore $(E^\ast,\tau^\ast)$ is a topological semilattice which
contains $E$ as a dense subsemilattice, a contradiction. The
obtained contradiction implies the assertion of the theorem.
\end{proof}

\begin{proposition}\label{proposition3.2}
Every linearly ordered topological pospace admits a structure of a
topological semilattice.
\end{proposition}

\begin{proof}
Let $X$ be a linearly ordered topological pospace and let
$\leqslant$ be a linear order on $X$. We define the semilattice
operation ``$*$'' on $X$ as follows:
\begin{equation*}
x*y=y*x=x \qquad \mbox{if} \quad x\leqslant y \quad (x, y\in X).
\end{equation*}
Since $\leqslant$ is a linear order on $X$, ``$*$'' is a
semilattice operation on $X$.

We observe that if $A\subseteq X$, then $A*A=A$. Let $x, y\in X$
be such that $x<y$ holds. Then $x*y=x$, and since $X$ is a
topological pospace, there exist neighbourhoods $U(x)$ and $U(y)$
of the points $x$ and $y$ in $X$, respectively, such that
\begin{equation*}
U(x)={\downarrow}U(x), \quad U(y)={\uparrow}U(y), \quad \mbox{and}
\quad U(x)\cap U(y)=\varnothing.
\end{equation*}
Let $V(x)$ be an arbitrary open neighbourhood of the point $x$ in
$X$. We put $W(x)=V(x)\cap U(x)$. Since $X$ is a linearly ordered
set we have
\begin{equation*}
U(y)*W(x)=W(x)*U(y)=W(x)\subseteq V(x).
\end{equation*}
Therefore $(X,*)$ is a topological semilattice.
\end{proof}

Proposition~\ref{proposition3.2} implies

\begin{proposition}\label{proposition3.3}
A linearly ordered topological  semilattice is $H$-closed if and
only if it is $H$-closed as a topological pospace.
\end{proposition}

A linearly ordered topological semilattice $E$ is called
\emph{complete} if every non-empty subset of $S$ has $\inf$ and
$\sup$.

In~\cite{GutikRepovs} Gutik and Repov\v{s} proved the following
theorem:

\begin{theorem}\label{Gutik-Repovs_theorem}
A linearly ordered topological semilattice $E$ is $H$-closed if
and only if the following conditions hold:
\begin{itemize}
    \item[$(i)$] $E$ is complete;
    \item[$(ii)$] $x=\sup A$ for $A={\downarrow}A\setminus\{ x\}$
    implies $x\in\operatorname{cl}_EA$, whenever
    $A\neq\varnothing$; and
    \item[$(iii)$] $x=\inf B$ for $B={\uparrow}B\setminus\{ x\}$
    implies $x\in\operatorname{cl}_EB$, whenever
    $B\neq\varnothing$
\end{itemize}
\end{theorem}

Propositions~\ref{proposition3.2}, \ref{proposition3.3} and
Theorem~\ref{Gutik-Repovs_theorem} imply the following:

\begin{corollary}\label{corollary3.4}
A linearly ordered topological pospace $X$ is $H$-closed if and
only if the following conditions hold:
\begin{itemize}
    \item[$(i)$] $X$ is a complete semilattice with the respect
    to the partial order on $X$;
    \item[$(ii)$] $x=\sup A$ for $A={\downarrow}A\setminus\{ x\}$
    implies $x\in\operatorname{cl}_XA$, whenever
    $A\neq\varnothing$; and
    \item[$(iii)$] $x=\inf B$ for $B={\uparrow}B\setminus\{ x\}$
    implies $x\in\operatorname{cl}_XB$, whenever
    $B\neq\varnothing$
\end{itemize}
\end{corollary}

A semilattice $S$ is called \emph{algebraically closed} (or
\emph{absolutely maximal}) if $S$ is a closed subsemilattice in
any topological semilattice which contains $S$ as a
subsemilattice~\cite{Stepp1975}. In \cite{Stepp1975} Stepp proved
that a semilattice $S$ is algebraically closed if and only if any
chain in $S$ is finite. Example~\ref{example3.5} shows that the
similar statement does not hold for an $H$-closed topological
pospace.

\begin{example}\label{example3.5}
Let $Y$ be an Hausdorff topological space with isolated points $a$
and $b$. We put $X=Y\setminus\{ a,b\}$. On $Y$ we define a partial
order $\leqslant$ as follows:
\begin{itemize}
    \item[1)] $x\leqslant x$ for all $x\in Y$;
    \item[2)] $x\leqslant b$ for all $x\in Y$; and
    \item[3)] $a\leqslant x$ for all $x\in Y$.
\end{itemize}
Obviously, $(Y,\leqslant)$ is a topological pospace, and moreover
$(Y,\leqslant)$ is an $H$-closed topological pospace if and only
if $X$ is an $H$-closed topological space.
\end{example}

We observe that the following conditions hold:
\begin{itemize}
    \item[1.] The partial order $\leqslant$ on $Y$ admits the
    lattice structure on $Y$:
    \begin{itemize}
        \item[1)] $b\vee x=b$ and $b\wedge x=x$ for all $x\in Y$;
        \item[2)] $x\wedge y=a$ and $x\vee y=b$ for all $x,y\in Y$;
        and
        \item[3)] $a\vee x=x$ and $a\wedge x=a$ for all $x\in Y$.
    \end{itemize}
    Hence the topological pospace $(Y,\leqslant)$ is both
    up-directed and down-directed. Also we observe that the
    lattice operations $\vee$ and $\wedge$ are not continuous in
    $Y$.
    \item[2.] The topological pospace $(Y,\leqslant)$ is both a
    $C_d$- and a $C_i$-space.
    \item[3.] The topological pospace $(Y,\leqslant)$ is upper and
    lower regularly ordered.
\end{itemize}

A partially ordered set $A$ is called a \emph{tree} if
${\downarrow}a$ is a chain for any $a\in A$.
Example~\ref{example3.6} shows that there exists an algebraically
closed (and hence $H$-closed) topological semilattice
$\mathscr{A}(\tau)$ which is a tree but $\mathscr{A}(\tau)$ is not
an $H$-closed topological pospace.

\begin{example}\label{example3.6}
Let $X$ be a discrete infinite space of cardinality $\tau$ and let
$\mathscr{A}(\tau)$ be the one-point Alexandroff compactification
of $X$. We put $\{\alpha\}=\mathscr{A}(\tau)\setminus X$ and fix
$\beta\in X$. On $\mathscr{A}(\tau)$ we define a partial order
$\leqslant$ as follows:
\begin{itemize}
    \item[1)] $x\leqslant x$ for all $x\in \mathscr{A}(\tau)$;
    \item[2)] $\beta\leqslant x$ for all $x\in \mathscr{A}(\tau)$;
    and
    \item[3)] $x\leqslant\alpha$ for all $x\in \mathscr{A}(\tau)$.
\end{itemize}
The partial order $\leqslant$ induces a semilattice operation
`$*$' on $\mathscr{A}(\tau)$:
\begin{itemize}
    \item[1)] $x*x=x$ for all $x\in \mathscr{A}(\tau)$;
    \item[2)] $\beta*x=x*\beta=\beta$ for all $x\in \mathscr{A}(\tau)$;
    \item[3)] $\alpha*x=x*\alpha=x$ for all $x\in \mathscr{A}(\tau)$;
    and
    \item[4)] $x*y=y*x=\beta$ for all distinct $x,y\in X$.
\end{itemize}
Since $X$ is a discrete subspace of $\mathscr{A}(\tau)$, $X$ with
induced from $\mathscr{A}(\tau)$ the semilattice operation is a
topological semilattice. By~\cite[Theorem~9]{Stepp1975} $X$ is an
algebraically closed semilattice, and hence it is an $H$-closed
topological semilattice. But $X$ a dense subspace of
$\mathscr{A}(\tau)$ and hence $X$ is not an $H$-closed pospace.
\end{example}


\section{Linearly ordered $H$-closed topological pospaces}

Let $C$ be a maximal chain of a topological pospace $X$. Then
$C=\bigcap_{x\in C}({\downarrow}x\cup{\uparrow}x)$, and hence $C$
is a closed subspace of $X$. Therefore we get the following:

\begin{lemma}\label{lemma2.3}
Let $K$ be a linearly ordered subspace of a topological pospace
$X$. Then $\operatorname{cl}_X(K)$ is a linearly ordered subspace
of $X$.
\end{lemma}

Since the conditions $(i)$---$(iii)$ of
Corollary~\ref{corollary3.4} are preserved by continuous monotone
maps, we have the following:

\begin{theorem}\label{theorem2.7}
Any continuous monotone image of a linearly ordered $H$-closed
topological pospace into a topological pospace is an $H$-closed
topological pospace.
\end{theorem}

Also Proposition~\ref{proposition2.1} follows from
Corollary~\ref{corollary3.4}.

\begin{proposition}\label{proposition2.1} Let $(X,\tau_X)$ be an
$H$-closed pospace of a linearly ordered topological pospace
$(T,\tau_T)$. Then the set ${\uparrow}{x}\cap X$
(${\downarrow}{x}\cap X$) contains a minimal (maximal) element for
any $x\in T$.
\end{proposition}

A subset $L$ of a linearly ordered set $X$ is called a
\emph{$L$-chain} in $X$ if ${\uparrow}x\cap{\downarrow}y\subseteq
L$ for any $x, y\in L$, $x\leqslant y$.

\begin{theorem}\label{theorem2.4} Let $X$ be a linearly ordered
topological pospace and let $L$ be a subspace of $X$ such that $L$
is an $H$-closed topological pospace and any maximal
$X{\setminus}L$-chain in $X$ is an $H$-closed topological pospace.
Then $X$ is an $H$-closed topological pospace.
\end{theorem}

\begin{proof} Suppose to the contrary that the topological
pospace $X$ is not $H$-closed. Then by Lemma~\ref{lemma2.3} there
exists a linearly ordered topological pospace $Y$ which contains
$X$ as a non-closed subspace. Without loss of generality we can
assume that $X$ is a dense subspace of a linearly ordered
topological pospace $Y$.

Let $x\in Y\setminus X$. The assumptions of the theorem imply that
the set $X\setminus L$ is a disjoint union of maximal
$X{\setminus}L$-chains $L_{\alpha}$, $\alpha\in{\mathscr A}$,
which are $H$-closed topological pospaces. Therefore any open
neighbourhood of the point $x$ intersects infinitely many sets
$L_{\alpha}$, $\alpha\in{\mathscr A}$.

Since any maximal $X{\setminus} L$-chain in $X$ is an $H$-closed
topological pospace, one of the following conditions holds:
\begin{equation*}
   {\uparrow}x\cap L\neq\varnothing \qquad \text{or} \qquad
   {\downarrow}x\cap L\neq\varnothing.
\end{equation*}
We consider the case when the sets ${\uparrow}x\cap L$ and
${\downarrow}x\cap L$ are nonempty. The proofs in the other cases
are similar.

By Proposition~\ref{proposition2.1} the set ${\uparrow}x\cap L$
contains a minimal element $x_m$ and the set ${\downarrow}x\cap L$
contains a maximal element $x_M$. Then the sets ${\uparrow}x_m$
and ${\downarrow}x_M$ are closed in $Y$ and, obviously,
$L\subset{\downarrow}x_M\cup{\uparrow}x_m$. Let $U(x)$ be an open
neighbourhood of the point $x$ in $Y$. We put
\begin{equation*}
    V(x)=U(x)\setminus\left({\downarrow}x_M\cup{\uparrow}x_m\right).
\end{equation*}
Then $V(x)$ is an open neighbourhood of the point $x$ in $Y$ which
intersects at most one maximal $S{\setminus} L$-chain $L_\alpha$,
a contradiction. Therefore $X$ is an $H$-closed topological
pospace.
\end{proof}

\begin{corollary}\label{corollary2.5}
Let $X$ be a linearly ordered topological pospace and let $L$ be a
subspace of $X$ such that $L$ is a compact topological pospace and
any maximal $X{\setminus}L$-chain in $X$ is a compact topological
pospace. Then $X$ is an $H$-closed topological pospace.
\end{corollary}

\begin{example}\label{example2.6}
Let ${\mathbb N}$ be the set of all positive integers. Let $\{
x_n\}$ be an increasing sequence in ${\mathbb N}$. Put ${\mathbb
N}^\ast=\{ 0\}\cup\{\frac{1}{n}\mid n\in{\mathbb N}\}$ and let
$\leqslant$ be the usual order on ${\mathbb N}^\ast$. We put
$U_n(0)=\{ 0\}\cup\{\frac{1}{x_k}\mid k\geqslant n\}$,
$n\in{\mathbb N}$. A topology $\tau$ on ${\mathbb N}^\ast$ is
defined as follows:
\begin{enumerate}
\item[a)] any point $x\in\mathbb{N}^\ast\setminus\{ 0\}$ is
isolated in ${\mathbb N}^\ast$; and
\item[b)] ${\mathscr B}(0)=\{ U_n(0)\mid n\in{\mathbb N}\}$ is
the base of the
topology $\tau$ at the point $0\in{\mathbb N}^\ast$.
\end{enumerate}
It is easy to see that $({\mathbb N}^\ast,\leqslant,\tau)$ is a
countable linearly ordered $\sigma$-compact locally compact
metrizable topological pospace and if $x_{k+1}>x_k+1$ for every
$k\in{\mathbb N}$, then $({\mathbb N}^\ast,\leqslant,\tau)$ is a
non-compact topological pospace.
\end{example}

By Corollary~\ref{corollary2.5} $({\mathbb
N}^\ast,\leqslant,\tau)$ is an $H$-closed topological pospace.
Also $({\mathbb N}^\ast,\leqslant,\tau)$ is a normally ordered (or
monotone normal) topological pospace, i.e. for any closed subset
$A={\downarrow}A$ and $B={\uparrow}B$ in $X$ such that $A\cap
B=\varnothing$ there exist open subsets $U={\downarrow}U$ and
$V={\uparrow}V$ in $X$ such that $A\subseteq U$, $B\subseteq V$,
and $U\cap V=\varnothing$~\cite{Nachbin1965}. Therefore for any
disjunct closed subsets $A={\downarrow}A$ and $B={\uparrow}B$ in
$X$ the exists a continuous monotone function $f\colon
X\rightarrow [0,1]$ such that $f(A)=0$ and
$f(B)=1$~\cite{Nachbin1965}.

Example~\ref{example2.6} implies negative answers to the following
questions:

\begin{enumerate}
    \item[$(i)$] Is a closed subspace of an $H$-closed topological
    pospace $H$-closed?
    \item[$(ii)$] Does any locally compact topological pospace embed into
    a compact topological pospace?
    \item[$(iii)$] Has any locally compact topological pospace a
    subbasis of open decreasing and open increasing subsets?
\end{enumerate}

Example~\ref{example1.2} shows that there exists a countably
compact topological pospace, whose space is $H$-closed. This
example also shows that there exists a countably compact totally
disconnected scattered topological pospace which is not embeddable
into a locally compact topological pospace.

\begin{example}\label{example1.2}
Let $X=[0,\omega_1)$ with the order topology (see
\cite[Example~3.10.16]{Engelking1989}), and let $Y=\{
0\}\cup\{\frac{1}{n}\mid n=1,2,3,\ldots\}$ with the natural
topology. We put $S=X\times Y$ with the product topology $\tau_p$
and the partial order $\preccurlyeq$:
\begin{equation*}
(x_1,y_1)\preccurlyeq(x_2,y_2) \quad \mbox{if}\quad y_1>y_2 \quad
\mbox{or} \quad  y_1=y_2 \mbox{ and } x_2\leqslant x_1.
\end{equation*}
We extend the partial order $\preccurlyeq$ onto
$S^*=S\cup\{\alpha\}$, where $\alpha\notin S$, as follows:
$\alpha\preccurlyeq\alpha$ and $\alpha\preccurlyeq x$ for all
$x\in S$, and define a topology $\tau$ on $S^*$ as follows. The
bases of topologies $\tau$ and $\tau_p$ at the point $x\in S$
coincide and the family $\mathscr{B}(\alpha)=\{
U_\beta(\alpha)\mid \beta\in\omega_1\}$ is the base of the
topology $\tau$ at the point $\alpha\in S^*$, where
\begin{equation*}
    U_\beta(\alpha)=\{\alpha\}\cup\big([\beta,\omega_1)
    \times\{1/n\mid n=1,2,3,\ldots\}\big).
\end{equation*}
Obviously, $(S^*,\preccurlyeq,\tau)$ is a Hausdorff non-regular
topological pospace. Proposition~3.12.5~\cite{Engelking1989}
implies that $(S^*,\tau)$ is an $H$-closed topological space. By
Corollary~3.10.14~\cite{Engelking1989} and
Theorem~3.10.8~\cite{Engelking1989} the topological space
$(S^*,\tau)$ is countably compact. Since every point of
$(S^*,\tau)$ has a singleton component, the topological space
$(S^*,\tau)$ is totally disconnected.

Let $A$ be a closed subset of $(S^*,\preccurlyeq,\tau)$ such that
$A\neq\{\alpha\}$. Then there exists $x\in[0,\omega_1)$ such that
$\tilde{A}=A\cap ([0,x]\times Y)\neq\varnothing$. Since
$[0,x]\times Y$ is compactum, $\tilde{A}$ is a compact topological
pospace, and hence $\tilde{A}$ contains a maximal element of
$\tilde{A}$. Let $x_m$ be a maximal element of $\tilde{A}$.
Definition of the topology $\tau$ on $S^*$ implies that
${\uparrow}x_m$ is an open subset in $(S^*,\tau)$. Then
${\uparrow}x_m\cap\tilde{A}=x_m$ and hence $x_m$ is an isolated
point of the space $\tilde{A}$ with the induced topology from
$(S^*,\tau)$. Therefore every closed subset of $(S^*,\tau)$ has an
isolated point in itself and hence $(S^*,\tau)$ is a scattered
topological space.
\end{example}

\begin{remark}\label{remark3.4}
The topological pospace $({\mathbb N}^\ast,\leqslant,\tau)$ from
Example~\ref{example2.6} admits the structure of a topological
semilattice:
\begin{equation*}
ab=\min\{ a, b\},\quad\text{ for }\quad a, b\in{\mathbb N}^\ast.
\end{equation*}

Also the topological pospace $(S^*,\preccurlyeq,\tau)$ from
Example~\ref{example1.2} admits the continuous semilattice
operation
\begin{equation*}
    (x_1,y_1)\cdot(x_2,y_2)=(\max\{ x_1,x_2\}, \max\{ y_1,y_2\})
\quad
\mbox{and} \quad
    (x_1,y_1)\cdot\alpha=\alpha\cdot(x_1,y_1)=\alpha,
\end{equation*}
for $x_1,x_2\in X$ and $y_1, y_2\in Y$.
\end{remark}

The following example shows that there exists a countable
$H$-closed scattered totally disconnected topological pospace
which has a non-$H$-closed maximal chain.

\begin{example}\label{example1.2a}
Let $X=\{ 1,2,3,\ldots\}$ be the positive integers with the
discrete topology, and let $Y=\{ 0\}\cup\{\frac{1}{n}\mid
n=1,2,3,\ldots\}$ with the natural topology. We put $T=X\times Y$
with the product topology $\tau_T$ and the partial order
$\preccurlyeq$:
\begin{equation*}
(x_1,y_1)\preccurlyeq(x_2,y_2) \quad \mbox{if}\quad y_1>y_2 \quad
\mbox{or} \quad  y_1=y_2 \mbox{ and } x_2\leqslant x_1.
\end{equation*}
We extend the partial order $\preccurlyeq$ onto
$T^*=T\cup\{\alpha\}$, where $\alpha\notin T$, as follows:
$\alpha\preccurlyeq\alpha$ and $\alpha\preccurlyeq x$ for all
$x\in T$, and define a topology $\tau^*$ on $T^*$ as follows. The
bases of topologies $\tau^*$ and $\tau_T$ at the point $x\in T$
coincide and the family $\mathscr{B}(\alpha)=\{ U_k(\alpha)\mid
k\in\{1,2,3,\ldots\}\}$ is the base of the topology $\tau^*$ at
the point $\alpha\in T^*$, where
\begin{equation*}
    U_k(\alpha)=\{\alpha\}\cup\big(\{k, k+1, k+2,\ldots\}
    \times\{1/n\mid n=1,2,3,\ldots\}\big).
\end{equation*}
It is obvious that $(T^*,\preccurlyeq,\tau^*)$ is a Hausdorff
non-regular topological pospace. Proposition~3.12.5
\cite{Engelking1989} implies that $(T^*,\tau^*)$ is an $H$-closed
topological space. Since every point of $(T^*,\tau^*)$ has a
singleton component, the topological space $(T^*,\tau^*)$ is
totally disconnected. The proof of the fact that $(T^*,\tau^*)$ is
a scattered topological pospace is similar to the proof of the
scatteredness of the topological pospace $(S^*,\preccurlyeq,\tau)$
in Example~\ref{example1.2}.

We observe that the set $L=\left(X\times
\{0\}\right)\cup\{\alpha\}$ with the induced partial order from
the topological pospace  $(T^*,\preccurlyeq,\tau^*)$ is a maximal
chain in $T^*$. The topology $\tau^*$ induces the discrete
topology on $L$. Corollary~\ref{corollary3.4} implies that $L$ is
not an $H$-closed topological pospace.
\end{example}

Theorem~\ref{theorem2.8} gives sufficient conditions on a maximal
chain of an $H$-closed topological pospace to be $H$-closed. We
shall say that a chain $L$ of a partially ordered set $P$ has a
\emph{${\downarrow}{\cdot}\max$-property}
(\emph{${\uparrow}{\cdot}\min$-property}) in $P$ if for any $a\in
P$ such that ${\downarrow}a\cap L\neq\varnothing$
(${\uparrow}a\cap L\neq\varnothing$) the chain {${\downarrow}a\cap
L$ (${\uparrow}a\cap L$) has a maximal (minimal) element. If the
chain of a partially ordered set $P$ has
${\downarrow}{\cdot}\max$- and ${\uparrow}{\cdot}\min$-properties,
then we shall call that $L$ has a
\emph{${\updownarrow}{\cdot}\operatorname{m}$-property}.

Similarly to \cite{McCartan1971, Priestley1972} we shall say that
a topological pospace $X$ is a \emph{$CC_i$-space} (resp.
\emph{$CC_d$-space}) if whenewer a chain $F$ of $X$ is closed,
${\uparrow}F$ (resp. ${\downarrow}F$) is a closed subset in $X$.

\begin{theorem}\label{theorem2.8}
Let $X$ be an $H$-closed topological pospace. If $X$ satisfies the
following properties:
\begin{itemize}
    \item[($i$)] $X$ is regularly ordered;
    \item[($ii$)] $X$ is $CC_i$-space; and
    \item[($iii$)] $X$ is $CC_d$-space,
\end{itemize}
then every maximal chain in $X$ with
${\updownarrow}{\cdot}\operatorname{m}$-property is an $H$-closed
topological pospace.
\end{theorem}

\begin{proof}
Suppose to the contrary that there is a non-$H$-closed chain $L$
with ${\updownarrow}{\cdot}\operatorname{m}$-property in $X$. Then
by Corollary~\ref{corollary3.4} at least one of the following
conditions holds:
\begin{itemize}
    \item[I)] the set $L$ is not a complete semilattice with the
    induced partial order from $X$;

    \item[II)] there exist a non-empty subset $A$ in $L$ with
    $x=\inf A$ such that $A={\uparrow}A\setminus\{ x\}$ and
    $x\notin\operatorname{cl}_L(A)$;

    \item[III)] there exist a non-empty subset $B$ in $L$ with
    $y=\sup B$ such that $B={\downarrow}B\setminus\{ y\}$ and
    $y\notin\operatorname{cl}_L(B)$.
\end{itemize}

Suppose that condition I) holds. Since topological space with dual
order to $\leqslant$ is a topological pospace, without loss
generality we can assume that there exists a subset $S$ in $L$
which does not have $\sup$ in $L$. Then the set ${\downarrow}S\cap
L$ does not have $\sup$ in $L$ either. Hence the set
$I=L\setminus{\downarrow}S$ does not have $\inf$ in $L$. We
observe that the maximality of $L$ implies that there exist no
$\inf I$ and $\sup S$ in such that $\sup S\leqslant\inf I$. Also
we observe that properties $(ii)-(iii)$ of $X$ and
Corollaries~\ref{corollary1.1b}, \ref{corollary1.1c} imply that
$I\neq\varnothing$. Therefore without loss of generality we can
assume that $S={\downarrow}S\cap L$, $I={\uparrow}I\cap L$ and $L$
is the disjoint inion of $S$ and $I$.

Since the set $S$ does not have $\sup$ in $L$, for any $s\in S$
there exists $t>s$, and hence there exist open neighbourhoods
$U(s)$ and $U(t)$ of the points $s$ and $t$ respectively such that
\begin{equation*}
    U(s)={\downarrow}U(s), \quad U(t)={\uparrow}U(t), \quad
    \mbox{and} \quad U(s)\cap U(t)=\varnothing.
\end{equation*}
Therefore we conclude that $S$ is an open subset of $L$ and $I$ is
a closed subset of $L$. Similarly we get $I$ is open in $L$ and
$S$ is closed in $L$.

By Lemma~1~\cite{Ward1954} for any $a,b\in X$ such that
$a\nleqslant b$ there exist open neighbourhoods $U(a)$ and $U(b)$
of the points $a$ and $b$ respectively such that
\begin{equation*}
    U(a)={\uparrow}U(a), \quad U(b)={\downarrow}U(b), \quad
    \mbox{and} \quad U(a)\cap U(b)=\varnothing.
\end{equation*}

We define the families
\begin{equation*}
    \mathscr{S}=\{(U(a),U(b))\mid a\nleqslant b,\, a\in S,\, b\in
X\}
\quad \mbox{and} \quad
    \mathscr{I}=\{(U(b),U(a))\mid a\nleqslant b,\, b\in I,\, a\in
X\}
\end{equation*}
as follows: $U(a)$ is an open neighbourhood of $a$ and $U(b)$ is
an open neighbourhood of $b$ such that
\begin{equation*}
    U(a)={\uparrow}U(a), \quad U(b)={\downarrow}U(b), \quad
    \mbox{and} \quad U(a)\cap U(b)=\varnothing.
\end{equation*}

Let $x\notin X$. We extend the partial order $\leqslant$ from $X$ onto
$X^\star=X\cup\{ x\}$ as follows:
\begin{itemize}
    \item[1)] $x\leqslant x$;
    \item[2)] $x\leqslant a$ for each $a\in I$;
    \item[3)] $x\leqslant z$ for $z\in X\setminus L$ if and only
    if
    there exists $e\in I$ such that $e\leqslant z$;
    \item[4)] $b\leqslant x$ for each $b\in S$; and
    \item[5)] $z\leqslant x$ for $z\in X\setminus L$ if and only
    if
    there exists $e\in S$ such that $z\leqslant e$.
\end{itemize}

We denote the topology on $X$ by $\tau_X$. On $X^\star$ we
determine a topology $\tau^\star$ as follows. For any point $y\in
X$ the bases of topologies $\tau^\star$ and $\tau_X$ at the point
$y$ coincide. For any $y\in X$ by $\mathscr{B}(y)$ we denote the
base of topology $\tau_X$ at the point $y$. We put
\begin{equation*}
    \mathscr{S}(x)=\{\{ x\}\cup U(a)\mid
    (U(a),U(b))\in\mathscr{S},\, a\nleqslant b,\, a\in S,\, b\in
    X\},
\end{equation*}
\begin{equation*}
    \mathscr{I}(x)=\{\{ x\}\cup U(b)\mid
    (U(b),U(a))\in\mathscr{S},\, a\nleqslant b,\, b\in I,\, a\in
    X\},
\end{equation*}
\begin{equation*}
    \mathscr{P}(x)=\mathscr{S}(x)\cup \mathscr{I}(x), \qquad
\mbox{and}
\end{equation*}
\begin{equation*}
\mathscr{B}(x)=\{ U_1\cap\cdots\cap U_n\mid
U_1,\ldots,U_n\in\mathscr{P}(x),\; n=1,2,3,\ldots\}.
\end{equation*}
Obviously, the conditions (BP1)--(BP3) of~\cite{Engelking1989}
hold for the family $\{{\mathscr B}(y)\}_{y\in X^\star}$ and hence
\break $\{{\mathscr B}(y)\}_{y\in X^\star}$ is a base of a
topology $\tau^\star$ at the point $y\in X^\star$. Since the set
$I$ has not an $\inf$, for finite many pairs
$(U(b_1),U(a_1)),\ldots, (U(b_n),U(a_n))\in\mathscr{I}(x)$ we have
that the set $U(b_1)\cap\ldots\cap U(b_n)\cap I$ is infinite.
Similarly, since the set $S$ does not have a $\sup$, for finite
many pairs $(U(a_1),U(c_1)),\ldots,
(U(a_n),U(c_k))\in\mathscr{S}(x)$ we have that the set
$U(a_1)\cap\ldots\cap U(a_k)\cap S$ is infinite. Therefore $x$ is
a non-isolated point of the topological space
$(X^\star,\tau^\star)$.

Further we shall show that $(X^\star,\tau^\star,\leqslant)$ is a
topological pospace. We consider three cases:
\begin{itemize}
    \item[1)] $y\in{\uparrow}I$;
    \item[2)] $y\in{\downarrow}S$; and
    \item[3)] $y\in X\setminus({\uparrow}I\cup{\downarrow}S)$.
\end{itemize}

If $y\in{\uparrow}I$ we have $x\leqslant y$, and since the set $I$
does not have an $\inf$, there exist $a,b\in I$ such that
$x\leqslant b<a\leqslant y$. Let $(U(b),U(a))\in\mathscr{I}$. Then
$U_b(x)=\{ x\}\cup U(b)$ and $U(y)=U(a)$ are open neighbourhoods
of the points $x$ and $y$ respectively such that
\begin{equation*}
    U_b(x)={\downarrow}U_b(x), \quad U(y)={\uparrow}U(y), \quad
    \mbox{and} \quad U_b(x)\cap U(y)=\varnothing.
\end{equation*}

If $y\in{\downarrow}S$ then we have $y\leqslant x$, and since the
set $S$ does not have a $\sup$, there exist $a,b\in S$ such that
$y\leqslant b<a\leqslant y$. Let $(U(a),U(b))\in\mathscr{S}$. Then
$U_a(x)=\{ x\}\cup U(a)$ and $U(y)=U(b)$ are open neighbourhoods
of the points $x$ and $y$ respectively such  that
\begin{equation*}
    U_a(x)={\uparrow}U_a(x), \quad U(y)={\downarrow}U(y), \quad
\mbox{and} \quad U_a(x)\cap U(y)=\varnothing.
\end{equation*}

Suppose $y\in X\setminus({\uparrow}I\cup{\downarrow}S)$. Then
$x\nleqslant y$. In the case ${\downarrow}y\cap L=\varnothing$ we
put $s=\min L$, and since $s\nleqslant y$, there exist open
neighbourhoods $U(s)$ and $U(y)$ of the points $s$ and $y$
respectively such that $(U(s),U(y))\in\mathscr{S}$. Then
$U_s(x)=\{ x\}\cup U(s)$ and $U(y)$ are open neighbourhoods of the
points $x$ and $y$ respectively such  that
\begin{equation*}
    U_s(x)={\uparrow}U_s(x), \quad U(y)={\downarrow}U(y), \quad
\mbox{and} \quad U_s(x)\cap U(y)=\varnothing.
\end{equation*}
If ${\downarrow}y\cap L\neq\varnothing$ then since the chain $L$
has ${\updownarrow}{\cdot}\operatorname{m}$-property, we put
$t=\max({\downarrow}y\cap L)$. Then since the set $S$ does not
have a $\sup$, there exists $m\in S\setminus{\downarrow}y$ such
that $t\leqslant m$. Obviously, $m\nleqslant y$. Then there exist
open neighbourhoods $U(m)$ and $U(y)$ of the points $m$ and $y$
respectively such that $(U(m),U(y))\in\mathscr{S}$. Then
$U_m(x)=\{ x\}\cup U(m)$ and $U(y)$ are open neighbourhoods of the
points $x$ and $y$ respectively such that
\begin{equation*}
    U_m(x)={\uparrow}U_m(x), \quad U(y)={\downarrow}U(y), \quad
\mbox{and} \quad U_m(x)\cap U(y)=\varnothing.
\end{equation*}

Suppose $y\nleqslant x$. Then by the statements dual to the
previous statements we get that there exist open neighbourhoods
$U(x)$ and $U(y)$ of the points $x$ and $y$ respectively such that
\begin{equation*}
    U(x)={\downarrow}U(x), \quad U(y)={\uparrow}U(y), \quad
\mbox{and} \quad U(x)\cap U(y)=\varnothing.
\end{equation*}

Therefore we get that $(X^\star,\tau^\star,\leqslant)$ is a
topological pospace and $X$  is a dense subspace of
$(X^\star,\tau^\star,\leqslant)$. This contradicts the assumption
that $X$ is an $H$-closed pospace.

Suppose that the statement II) holds, i.~e., that there exists an
open neighbourhood $O(x)$ of $x=\inf A$ such that $O(x)\cap
A=\varnothing$. Since $X$ is a topological pospace, for every
$t\in L\setminus(A\cup\{ x\})$ there exists an open decreasing
neighbourhood $U(t)$ such that $x\notin U(t)$, and hence $A$ is a
closed subset of a topological space $X$.

Let $a\notin X$. We pu $X^\dagger=X\cup\{ a\}$ and extend the
partial order $\leqslant$ from $X$ onto $X^\dagger$ as follows:
\begin{itemize}
    \item[1)] $a\leqslant a$;
    \item[2)] $a\leqslant b$ for each $b\in A$;
    \item[3)] $a\leqslant z$ for $z\in X\setminus L$ if and only
    if there exists $e\in A$ such that $e\leqslant z$;
    \item[4)] $b\leqslant a$ for each $b\in L\setminus A$; and
    \item[5)] $z\leqslant a$ for $z\in X\setminus L$ if and only
    if there exists $e\in L\setminus A$ such that $z\leqslant e$.
\end{itemize}

By Lemma~1~\cite{Ward1954} we define the family
\begin{equation*}
    \mathscr{A}=\{(U(t),U(b))\mid b\nleqslant t,\, t\in A,\, b\in
X\}
\end{equation*}
as follows: $U(t)$ and $U(b)$ are open neighbourhoods of the
points $t$ and $b$ respectively such that
\begin{equation*}
    U(t)={\downarrow}U(t), \quad U(b)={\uparrow}U(b), \quad
\mbox{and} \quad U(t)\cap U(b)=\varnothing.
\end{equation*}

Since $X$ is a regularly ordered $CC_i$- and $CC_d$-space, we
determine the family
\begin{equation*}
    \mathscr{V}=\{(V_t(A),V(t))\mid t\notin
X\setminus{\uparrow}A\}
\end{equation*}
as follows: $V_t(A)$ and $V(t)$ are open neighbourhoods of the set
$A$ and the point $t$ respectively such that
\begin{equation*}
    V_t(A)={\uparrow}V_t(A), \quad V(t)={\downarrow}V(t), \quad
\mbox{and} \quad V_t(A)\cap V(t)=\varnothing.
\end{equation*}

We determine the topology $\tau^\dagger$ on $X^\dagger$ as
follows. Let $\tau_X$ be the topology on $X$. For any point $y\in
X$ the base of topologies $\tau^\dagger$ and $\tau_X$ at the point
$y$ coincide. For every $y\in X$ by $\mathscr{B}(x)$ we denote a
base of the topology $\tau_X$ at the point $y$. We put
\begin{equation*}
    \mathscr{A}(x)=\{\{a\}\cup U(t)\mid (U(t),U(b))\in\mathscr{A},
b\nleqslant t, \, t\in A,\, b\in X\},
\end{equation*}
\begin{equation*}
    \mathscr{V}(x)=\{\{ a\}\cup V(A)\mid
(V(A),V(t))\in\mathscr{V}, \, t\in X\setminus{\uparrow}A\},
\end{equation*}
\begin{equation*}
    \mathscr{P}(x)=\mathscr{A}(x)\cup \mathscr{V}(x), \qquad
    \mbox{and}
\end{equation*}
\begin{equation*}
\mathscr{B}(x)=\{ U_1\cap\cdots\cap U_n\mid
U_1,\ldots,U_n\in\mathscr{P}(x),\; n=1,2,3,\ldots\}.
\end{equation*}
Obviously, the conditions (BP1)--(BP3) of~\cite{Engelking1989}
hold for the family $\{{\mathscr B}(y)\}_{y\in X^\dagger}$ and
hence \break $\{{\mathscr B}(y)\}_{y\in X^\dagger}$ is a base of a
topology $\tau^\dagger$ at the point $y\in X^\dagger$. Since for
finitely many pairs
$(U(t_1),U(b_1)),\ldots,(U(t_k),U(b_k))\in\mathscr{A}$ the
intersection $U(t_1)\cap\ldots\cap(U(t_k)\cap A$ is an infinite
set, $a$ is a non-isolated point of the topological pospace
$(X^\dagger,\tau^\dagger)$.

Further we shall show that $(X^\dagger,\tau^\dagger,\leqslant)$ is
a topological pospace. We consider three cases:
\begin{itemize}
    \item[1)] $y\in{\uparrow}A$;
    \item[2)] $y\in{\downarrow}(L\setminus A)$;
    \item[3)] $y\in X\setminus({\uparrow}A\cup{\downarrow}(L\setminus
    A)$.
\end{itemize}

If $y\in{\uparrow}A$ then we have $a\leqslant y$, and since the
set $A$ does not contain $\min$, there exist $c,d\in A$ such that
$a<c<d\leqslant y$. Let $(U(c),U(d))\in\mathscr{A}$. Then
$U_c(a)=\{ a\}\cup U(c)$ and $U(y)=U(d)$ are open neighbourhoods
of the point $a$ and $y$ respectively such that
\begin{equation*}
    U_c(a)={\downarrow}U_c(a), \quad U(y)={\uparrow}U(y), \quad
\mbox{and} \quad U_c(a)\cap U(y)=\varnothing.
\end{equation*}

If $y\in{\downarrow}(L\setminus A)$ then we have $y\leqslant x<a$.
Let $(V_x(A),V(x))\in\mathscr{V}$. Then $V_x(a)=\{ a\}\cup V_x(A)$
and $V(y)=V(x)$ are open neighbourhoods of the points $a$ and $y$
respectively such that
\begin{equation*}
    V_x(a)={\uparrow}V_x(a), \quad V(y)={\downarrow}V(y), \quad
\mbox{and} \quad V_x(a)\cap V(y)=\varnothing.
\end{equation*}

Suppose $y\in X\setminus({\uparrow}A\cup{\downarrow}(L\setminus
A)$. Then $a\nleqslant y$. Let $(V_y(A),V(y))\in\mathscr{V}$. Then
$V_y(a)=\{ a\}\cup V_y(A)$ and $V(y)$ are open neighbourhoods of
the points $a$ and $y$ respectively such that
\begin{equation*}
    V_y(a)={\uparrow}V_y(a), \quad V(y)={\downarrow}V(y), \quad
\mbox{and} \quad V_y(a)\cap V(y)=\varnothing.
\end{equation*}

Suppose $y\nleqslant a$. In the case ${\uparrow}y\cap
A=\varnothing$ since by Corollary~\ref{corollary1.1c} $\max A$
exists, we put $s=\max A$. Then $y\nleqslant s$ and hence there
exists $(U(s),U(y))\in\mathscr{A}$. Then $U_s(a)=\{ a\}\cup U(s)$
and $U(y)$ are open neighbourhoods of the points $a$ and $y$
respectively such that
\begin{equation*}
    U_s(a)={\downarrow}U_s(a), \quad U(y)={\uparrow}U(y), \quad
\mbox{and} \quad V_s(a)\cap V(y)=\varnothing.
\end{equation*}

If ${\uparrow}y\cap A\neq\varnothing$ then since the chain $L$ has
${\updownarrow}{\cdot}\operatorname{m}$-property, we put
$t=\min({\uparrow}y\cap L)$. Then since the set $A$ does not
contain a $\min$, there exists $m\in A\setminus{\uparrow}y$ such
that $m\leqslant t$. Obviously $y\nleqslant m$, and hence there
exists $(U(m),U(y))\in\mathscr{A}$. Then $U_m(a)=\{ a\}\cup U(m)$
and $U(y)$ are open neighbourhoods of the points $a$ and $y$
respectively such that
\begin{equation*}
    U_m(a)={\downarrow}U_m(a), \quad U(y)={\uparrow}U(y), \quad
\mbox{and} \quad V_m(a)\cap V(y)=\varnothing.
\end{equation*}

Therefore we get that $(X^\dagger,\tau^\dagger,\leqslant)$ is a
topological pospace and $X$  is a dense subspace of
$(X^\dagger,\tau^\dagger,\leqslant)$. This contradicts the
assumption that $X$ is an $H$-closed pospace.

In case III) we get a similar contradiction to II).

The obtained contradictions imply the statement of the theorem.
\end{proof}

\begin{remark}\label{remark3.4a}
We observe that the topological pospace
$(T^*,\preccurlyeq,\tau^*)$ from Example~\ref{example1.2a} is not
regularly ordered and is not a $CC_i$-space. Also the topological
pospace $(T^*,\preccurlyeq,\tau^*)$ admits the continuous
semilattice operation
\begin{equation*}
    (x_1,y_1)\cdot(x_2,y_2)=(\max\{ x_1,x_2\}, \max\{ y_1,y_2\})
\quad
\mbox{and} \quad
    (x_1,y_1)\cdot\alpha=\alpha\cdot(x_1,y_1)=\alpha,
\end{equation*}
for $x_1,x_2\in X$ and $y_1, y_2\in Y$. Therefore a maximal chain
of an $H$-closed topological semilattice is not necessarily to be
an $H$-closed topological semilattice.
\end{remark}

\section*{Acknowledgements}

This research was support by the Slovenian Research Agency grants
P1-0292-0101-04, \break J1-9643-0101 and BI-UA/07-08/001.


\end{document}